 \documentclass[11pt]{amsart}
\usepackage{amssymb}
\newtheorem{theorem}{Theorem}
\newtheorem{lemma}[theorem]{Lemma}

\newtheorem{corollary}[theorem]{Corollary}

\theoremstyle{remark}
\newtheorem{remark}{Remark}
\newtheorem{example}{Example}

\newcommand{\calQ}{\mathcal{Q}}

\newcommand{\bbP}{\mathbb{P}}

\newcommand{\bbZ}{\mathbb{Z}}

\newcommand{\bfF}{\mathbf{F}}

\newcommand{\la}{\langle}
\newcommand{\ra}{\rangle}

\newcommand{\GL}{\textup{GL}}

\newcommand{\Pic}{\textup{Pic}}

\newcommand{\Cr}{\textup{Cr}}
\newcommand{\Aut}{\textup{Aut}}

\begin{document}
 \title[On elements of order $p^n$]{On  elements of order $p^s$  in the plane Cremona group  over a  field of characteristic $p$}

\author{Igor V. Dolgachev}

\address{Department of Mathematics, University of Michigan, 525 E. University Av., Ann Arbor, Mi, 49109, USA}
\email{idolga@umich.edu}

\dedicatory{To the memory of Vasily Iskovskikh}

\begin{abstract}
We show that the plane Cremona group over a  field of characteristic $p > 0$ does not contain elements of  order of power of $p$ larger than 2 and it does not contain elements of order $p^2$ unless $p =2$. Also we describe  conjugacy classes of elements of order $4$.
\end{abstract}

\maketitle

\bigskip\noindent
\section{Introduction}
The classification of conjugacy classes of elements of finite order
$\ell$ in the plane Cremona group $\Cr_2(k)$ over an algebraically
closed field $k$ of characteristic 0 has been known for more than a
century. The possible orders of elements not conjugate to a
projective transformation   are $2,3,4,5,6,8,9,10,$
$12,14,15,18,20,24, 30$ and any even order is realized by a de Jonqui\`eres transformation (see \cite{DI} and
historic references there). Much less is known in the case when  
$k$ is of positive characteristic $p$ and the order is divisible by  $p$. 

In this note we prove the following Main Theorem.

\begin{theorem}\label{main} Let $k$ be a field of characteristic $p > 0$. Then the group $\Cr_2(k)$ does not contain elements of order $p^s$ with $s > 2$.
\end{theorem}

We will also describe conjugacy classes of elements of order $p^2$ over algebraically closed field of characteristic $p > 0$.

\medskip

I thank J.-P. Serre for asking about the existence of elements of order 8 in $\Cr_2(k)$ over a field of characteristic  $2$ and his numerous comments on the previous versions of the paper. The question had initiated the present paper.   

For more than 45 years, Vasya Iskovskikh had been a friend, a collaborator on several papers and an inspiring guide in the area of birational geometry. He will be greatly missed.

\section{Conic bundles} It is clear that in the proof of Main Theorem,  we may assume that $k$ is an algebraically closed field of characteristic $p > 0$.  On several occasions I refer to \cite{DI} where the ground field was assumed to be the field of complex numbers. The proofs of the facts which I will use 
extend to our case.

Let $\sigma\in \Cr_2(k)$ be of order $p^s$. A standard argument  (see \cite{DI}) shows that $\sigma$ acts biregularly on one of the following rational surfaces $X$ 

\begin{itemize}
\item[(i)] $X$ has  a structure of a conic bundle $f:X\to \bbP_k^1$ with $m\ge 0$ singular fibres,
\item[(ii)] $X$ is a Del Pezzo surface of degree $d$.
\end{itemize}

Moreover, we may assume that $X$ is $\sigma$-minimal, i.e. $\Pic(X)^\sigma$ is of rank 2 in the first case and of rank 1 in the second case. This is equivalent to that any  $\sigma$-equivariant birational morphism $X\to X'$ must be an isomorphism. When $X$ is $\sigma$-minimal,  we say that $\sigma$ acts minimally on $X$. 

We start from the first case. Recall the following well-known fact.

\begin{lemma}\label{jordan} Let $\sigma$ be an element of order $p^s$ in $\Aut(\bbP_k^r)$. Then $s<1+\log_p(r+1)$.
\end{lemma}

\begin{proof} Let $A\in \GL_{r+1}(k)$ represent $\sigma$ and $A^{p^s} = cI_{r+1}$ for some constant $c$. Multiplying $A$ by $c^{\frac{1}{p^s}}$ we may assume that $A^{p^s} = I_{r+1}$ but $A^{p^{s-1}} \ne I_{r+1}$. Since  $k^*$ does not contain non-trivial $p$-th roots of unity, we can reduce $A$ to the Jordan  form with 1 at the diagonal.  Obviously $A^{p^{s-1}} = I_{r+1}+(A-I_{r+1})^{p^{s-1}}$. Since, for any Jordan block-matrix $J$  with zeros at the diagonal we have $J^{r+1} = 0$, we get  $p^{s-1} <  r+1$. The assertion follows. 
\end{proof}

\begin{corollary} Let $f:X\to \bbP_k^1$ be a conic bundle and $\sigma$ be an automorphism of $X$ of order $p^s$ preserving the conic bundle.   Then $s\le 2$. 
\end{corollary}

\begin{proof} Let $\bar{g}$ be the image of $\sigma$ in the automorphism group of the base of the fibration. By the previous lemma $\bar{\sigma}^p = 1$. Thus $\sigma^p$ acts identically on the base and hence acts on the general fibre of $f$. By Tsen's Theorem, the latter is isomorphic to the projective line over the function field of the base. Applying the lemma again we obtain that $\sigma^{p^2} = 1$.
\end{proof}

This checks the theorem in the case of a conic bundle. Let us give a closer look at elements of order $p^2$.

\begin{theorem} Let $\sigma$ be a minimal automorphism of order $p^2$ of a conic bundle $X\to \bbP_k^1$. Then $p = 2$.
\end{theorem}

\begin{proof}  Let $m = K_X^2-8$ be the number of  singular fibres of the conic bundle. Assume first that $m = 0$, i.e. $\pi:X\to \bbP_k^1$ is a minimal ruled surface $\bfF_n$. If $n = 1$, the surface is not $\sigma$-minimal. If $n = 0$, the automorphism group of $\bfF_0\cong \bbP_k^1\times \bbP_k^1$ preserving one of the rulings is isomorphic to $\Aut(\bbP_k^1)\times \Aut(\bbP_k^1)$. It does not contain elements of order $p^2$. 

So we may assume that $n \ge 2$. The automorphism group $\Aut(X)$ of the surface $\bfF_n$ is well-known (see \cite{DI}, \S 4.4). By blowing down the exceptional section, we obtain that $\Aut(X)$ is isomorphic to the group of automorphisms of the weighted projective plane $\bbP(1,1,n)$ with coordinates $t_0,t_1$ of degree $1$ and coordinate $t_2$ of degree $n$. Any automorphism $g$ of $\bbP(1,1,n)$ can be given by the formula
$$\sigma:(t_0,t_1,t_2)\mapsto (at_0+bt_1,ct_0+dt_1,et_2+f_n(t_0,t_1)),$$
where $f_n$ is a binary form of degree $n$. In our case we can  change the  coordinates to assume that $a=  b = d = 1, c  = 0$. By iterating, we get $e^{p^s} = 1$, hence $e = 1$. Also 
$$\sigma^p:(t_0,t_1,t_2) = (t_0,t_1,t_2+\sum_{j=0}^{p-1}f_n(t_0+jt_1,t_1).$$
Let $\bar{\sigma}$ be the transformation $(t_0,t_1)\mapsto (t_0+t_1,t_1)$. Since 
$\sum_{i=0}^{p-1}\bar{\sigma}^i = 0$, we get that the sum in above is equal to zero, hence $\sigma^p = 1$.
Thus there are no automorphisms of order $p^2$.

Assume now that $m > 0$, i.e. $X$ is obtained from a minimal ruled surface $\bfF_n$ by blowing up $m$ points. If $n > 0$,  the proper transform of the exceptional section  of $\bfF_n$ is a section of the conic bundle with negative self-intersection. If $n = 0$, the proper transform of a section of $\bfF_0$ passing through a point we blow up, is a section  with negative self-intersection. So, in any case we have a section of the conic bundle with negative self-intersection. It intersects a component of a singular fibre at its nonsingular point. Since $X$ is $\sigma$-minimal, $\sigma$ cannot fix this component, so $\sigma(E) \ne E$. By Lemma \ref{jordan},   $\sigma^p$ acts identically on the base of the conic bundle. Since $p > 2$,  $\sigma^p$ cannot switch components of singular fibres,  hence it must act identically on $\Pic(X)$. Since an irreducible  curve with negative self-intersection does not move in a linear system, $\sigma^p$ fixes $E$ and $\sigma(E)$. But in characteristic $p > 0$ an automorphism of order $p$ of a general fibre has only one fixed point.  This shows that $\sigma^p = 1$  if $p > 2$. 
\end{proof}

\begin{example} Recall that $\Cr_2(k)$ contains a subgroup of de Jonqui\`eres transformations of the form $(x,y) \mapsto \bigl(\frac{\alpha x+\beta}{\gamma x+\delta},\frac{a(x)y+b(x)}{c(x)y+d(x)}\bigr)$. Each element of finite order in this subgroup is realized as an automorphism of a conic bundle. Assume $p = 2$. Without loss of generality we may assume that $x\mapsto x+1$. 

Let $a(x) = d(x) =  xP(x),$ where $P(x$ is a polynomial of degree $n$ without multiple roots.  Let 
$b(x) = P(x)P(x+1), c(x) = x(x+1).$ We have  
$$a(x)a(x+1)+b(x)c(x+1) = a(x)a(x+1)+b(x+1)c(x) = 0.$$  With this choice, we have 
$\sigma^2:(x,y)\mapsto (x, R(x)/y),$ where 
$$R(x) = \tfrac{a(x+1)b(x)+a(x)b(x+1)}{a(x)c(x+1)+a(x+1)c(x)} = 
 \tfrac{P(x)P(x+1)}{x(x+1)}.$$
Replacing $y$ with $x(x+1)y$, we obtain the de Jonqui\`eres involution $(x,y) \mapsto (x,P(x)P(x+1)/y)$. It is known that it  is realized as a minimal automorphism of a conic bundle with the number $m$ of singular fibres equal to the degree of $P(x)P(x+1)$. On the other hand, it is known that for $m\ge 8$ a  minimal automorphism of such a conic bundle is not conjugate to neither a projective automorphism, nor a minimal automorphism of a Del Pezzo surface, nor a minimal automorphism of a conic bundle with number of singular fibres different from $m$ (see Corollary 7.11 in \cite{DI}). Thus we  have constructed a countable set of  conjugacy classes of elements of order 4 in $Cr_2(k)$.
\end{example}

\section{Del Pezzo surfaces of degree $\ge 3$}
Now we consider the case when $\sigma$ is an automorphism of order $p^s$ of a Del Pezzo surface $X$ of degree $d:= K_X^2 \ge 4$. 

If $d = 9$, $X = \bbP_k^2$ and by Lemma \ref{jordan} we get $s\le 2$.  All elements of order $p^2$ are conjugate in $\Aut(\bbP_k^2)$. 

If $d = 8$, then $X \cong \bbP_k^1\times \bbP_k^1$ because the ruled surface $\bfF_1$ is not $\sigma$-minimal. We know that $\Aut(\bfF_0)$ contains a subgroup of index 2 isomorphic to $\Aut(\bbP_k^1)\times \Aut(\bbP_k^1)$. Applying Lemma \ref{jordan} we obtain $s =1$ if $p \ne 2$, and $s\le 2$ otherwise.  The automorphism  of $X$ given in affine coordinates by $(x,y)\mapsto (y+1,x)$ is of order 4. 

If $d = 7$, the surface is not $\sigma$-minimal since it is obtained by blowing up two points in $\bbP_k^2$, the proper transform of the line joining the points is a $\sigma$-invariant $(-1)$-curve.

Assume $d = 6$. Then $\Aut(X)$ is isomorphic to the semi-direct product $T\rtimes G$, where $T \cong k^*{}^2$ is a 2-dimensional torus and $G$ is a dihedral group $D_{12} \cong (\bbZ/2\bbZ)\times S_3$. Since $T$ does not contain elements of order $p$ and $D_{12}$ does not contain elements of order $p^s, s > 1$, we obtain that  the only possibility is $s = 1$ and $ p = 2,3$.

Assume $d = 5$. It is known that   $\Aut(X)$  acts faithfully on the Picard group of $X$ of a Del Pezzo surface of degree $\le 5$. Via this action it becomes isomorphic to a subgroup of the Weyl group $W(A_4) \cong S_5$. Thus $s=1$ unless $p = 2$ and $s = 2$. The group $W(A_4)$ acts on $K_X^\perp \cong \bbZ^4$ 
via its standard irreducible representation on $\{(a_1,\ldots,a_5)\in \bbZ^5:a_1+\ldots+a_5 = 0\}$. A cyclic permutation of order 4 has a fixed vector. This shows that $X$ is not $\sigma$-minimal.

Assume $d = 4$. In this case  $\Aut(X)$ is  isomorphic to a subgroup of the Weyl group $W(D_5) \cong (\bbZ/2\bbZ)^4\rtimes S_5$. Thus an automorphism of order $p^s$ with $s> 1$ may exist only if $p = 2$.

 It is known that $X$ is  isomorphic to the blow-up of 5 points $p_1,\ldots,p_5$ in the plane, no three among them are collinear. The surface admits 5 pairs $(|C_i|, |C_i'|)$ of pencils of conics  in the anti-canonical embedding $X\hookrightarrow \bbP_k^4$. The pencil $|C_i|$ is the proper transform of the pencil of lines through the point $p_i$ and the pencil $|C_i'|$ is the proper transform of the pencil of conics through the points $p_j, j\ne i$. Since $C_i+C_i' \sim -K_X$,  the 
Weyl group permutes the 5 pairs of the divisor classes $[C_i], [C_i]'$ and switches $[C_i]$ with $[C_i']$ in even pairs of them (see \cite{DI}, Proposition 6.6). It is known that the anti-canonical linear system $|-K_X|$ maps $X$ isomorphically onto the intersection of two quadrics in $\bbP_k^4$. Under the multiplication map $|C_i|\times  |C_i'| \to |-K_X|$, the two pencils  generate a hyperplane $H_i$ in $|-K_X|$ and the map $f_i\times f_i':X\to \bbP_k^1\times \bbP_k^1$ defined by the two pencils is equal to the composition of the anti-canonical map and the projection from the point $h_i\in |-K_X|^*$ corresponding to the hyperplane $H_i$. Since the image of $X$ under this projection is a nonsingular quadric, we see that the center of the projection lies on a singular quadric $Q_i$ of corank 1 in  the pencil $\calQ$ of quadrics containing $X$. Conversely, every such quadric defines a degree 2 map $f:X \to \bbP_k^1\times \bbP_k^1$, and the pre-images of the ruling define a pair of pencils of conics on $X$. Thus we see that the pencil of quadrics $\calQ$ contains exactly five  singular quadrics. Any automorphism $\sigma$ of $X$ acts on the pencil $\calQ$ leaving the set of five quadrics invariant.  Its square $\sigma^2$ acts identically on the pencil and hence leaves invariant all pairs of conic pencils. Since the divisor classes $[C_i]$ together with $K_X$ generate $\Pic(X)$, we obtain that $\sigma^4$ acts identically on $\Pic(X)$, hence it is the identity.

\begin{remark} Another proof of non-existence of an automorphism of order 8 on a Del Pezzo surface of degree $4$ was suggested by J.-P. Serre. It is known that an element of order 8 in $W(D_5)$ has trace equal to $-1$ in the root lattice. Since the latter is isomorphic to $K_X^\perp$, the automorphism of order 8 has trace 0 in $\Pic(X)$ and hence in the second cohomology group with $\ell$-adic coefficients. Thus the Lefschetz number of $\sigma$ is equal to 2, and hence, by the Lefschetz-fixed-point formula, $\sigma$ has  a fixed point. Blowing it up we get an automorphism of order 8 of a cubic surface.  Since any automorphism of a cubic surface is the restriction of an automorphism of $\bbP_k^3$, applying Lemma \ref{jordan} we find a contradiction.
\end{remark}

Let us summarize what we have learnt.

\begin{theorem} A Del Pezzo surface of degree $\ge 4$ does not contain elements of order $p^3$. An automorphism  of order $p^2$ not conjugate to a projective automorphism in $\Cr_2(k)$ exists only if $p = 2$. It is minimally realized on $X = \bbP_k^1\times \bbP_k^1$ or on   a Del Pezzo surface of degree  4.
\end{theorem}

Note that any automorphism of order 4 of $\bbP_k^1\times \bbP_k^1$ has a fixed point, and the projection from this fixed point makes it conjugate to a projective transformation.

Assume now that $d = 3$, i.e. $X$ is a cubic surface embedded in  $\bbP_k^3$ by the anti-canonical linear system $|-K_X|$. In this case $\Aut(X)$ is isomorphic to a subgroup of the Weyl group $W(E_6)$ of a simple root lattice of type $E_6$. By Corollary 6.11 from \cite{DI}, all elements of order $p^s, s > 1$, in $W(E_6)$ have an invariant vector in the lattice $E_6 \cong K_X^\perp$ unless $p^s = 9$. Thus we have to consider the existence of an automorphism $\sigma$ of order 9 of a cubic surface over a field of characteristic $p = 3$.  

The following argument was suggested to me by J.-P. Serre. It follows from the classification of conjugacy classes of elements of $W(E_6)$ that the trace of $\sigma$ in its action in $K_X^\perp$ is equal to 0. Thus the Lefschetz number of $\sigma$ in in the $\ell$-adic cohomology of $X$ is equal to 3. This implies that $\sigma$ has a fixed point $x_0$. Since  $\sigma$ acts trivially on $|-K_X-x_0| \cong \bbP_k^2$, we obtain that it acts trivially on $|-K_X| \cong \bbP_k^3$.  

We have proved the following.

\begin{theorem} A cubic surface does not admit minimal automorphisms of order $p^s$ with $s > 1$.
\end{theorem}

\section{Del Pezzo surfaces of degree 2} 

It is known (see \cite{Dem}) that the linear system $|-K_X|$ defines a degree 2 map $f:X\to \bbP_k^2$. The map must be finite since $-K_X$ is ample. It is also a separable map because otherwise $X$ must be homeomorphic to $\bbP_k^2$, but comparing the $l$-adic Betti numbers we find this impossible. The cover $f$ is a Galois cover  with order 2 cyclic Galois group $\la \gamma\ra$. The automorphism $\gamma$ of $X$ is called the \emph{Geiser involution}. For any divisor $D$ we have 
$$D+\gamma^*(D)\sim (D\cdot K_X)K_X.$$
This implies that $\gamma^*$ acts on $K_X^\perp$ as the minus identity.  The lattice $K_X^\perp$ is isomorphic to the root lattice of type $E_7$, and the isometry $\gamma^*$  generates  the center of the Weyl group $W(E_7)$.

 It follows from the classification of conjugacy classes in $W(E_7)$ that  for any automorphism of order $p^s, s > 1$,  the rank of $\Pic(X)^\sigma$ is greater than 1, unless $p = s = 2$. So, it suffices to consider the latter case. All such automorphisms form one conjugacy class (of type $2A_3+A_1$ in the notation from \cite{DI}). It follows from the description of degree 2 covers of smooth varieties (see \cite{CD}, Chapter 0) that $X$ is isomorphic to a surface $\bbP(1,1,1,2)$ given by an equation
$$u^2+a_2(x,y,z)u+a_4(x,y,z) = 0,$$
where $a_2,a_4$ are homogeneous forms of degree 2 and 4. Since the anticanonical map is separable we have $a_2 \ne 0$. An automorphism $\sigma$ of order $4$ acts linearly in $\bbP_k^2 = |-K_X|^*$ leaving the branch curve $V(a_2) $ invariant. If $V(a_2)$ is an irreducible conic, then $\sigma^2$ is identical on the conic, and hence it is identical on $\bbP_k^2$. This implies that $\sigma^2$ is the Geiser  involution $u\mapsto u+a_2$. However,  the Weyl group $W(E_7)$ does not contain square roots of the Geiser involution. Suppose now that $V(a_2)$ is reducible. If it is not a double line, we can choose projective coordinates $x,y,z$ to assume that $a_2 = xy$. Then $\sigma^2$ must change $z$ to $z+ax+by$ and leave $x,y$ unchanged. This forces $a_4$ to be invariant with respect to this transformation. Writing 
$$a_4 = l_0z^4+z^3l_1+z^2l_2+zl_3+l_4,$$ where $l_i$ are binary forms in $x,y$, we find that $l_1 = 0$. This implies that the point $(x,y,z,u) = (0,0,1,0)$ is a singular point on the surface. Thus $\sigma^2$ must be the Geiser involution and we finish as in the previous case. Finally we may assume that the equation of $X$ looks like
$u^2+x^2u+a_4 = 0$. In this case, $\sigma^*(x) = x$ and we may assume that $\sigma$ acts on the variables $x,y,z$ by $x\mapsto x, y\mapsto y+x, z\mapsto z+y$. The  polynomial $a_4(x,y,z)$ must be invariant with respect to the coordinate change $\sigma^2:(x,y,z)\mapsto (x,y,z+x)$. It is easy to see that 
the ring of polynomials in $x,z$ invariant with respect $(x,z)\mapsto (x,z+x)$ is generated by $x$ and $z(z+x)$.   This implies that $a_4$ can be written as a polynomial in $z(z+x),x,y$
$$a_4 =  cz^2(z+x)^2+z(z+x)g(x,y)+h(x,y).$$
It is immediate to check that the point $(x,y,z,u) = (0,0,1,\sqrt{c})$ is a singular point of the surface. 

To sum up,  a Del Pezzo surface of degree 2 does not contain minimal automorphisms of order $p^s, s > 1$.

\begin{remark} Another argument to prove that a Del Pezzo surface $X$ of degree 2 has no elements of order 8 was suggested by J.-P. Serre.  We use that $W(E_7) = W(E_7)^+\times \la w_0\ra$, where $w_0$ generates the center of $W(E_7)$. In the faithful representation $\rho:\Aut(X)\to W(E_7)$, the image of the Geiser involution $\gamma$ is equal to $w_0$. This implies that a subgroup $G$ of order 8 of $\Aut(X)$ is isomorphic to a subgroup of $A\times \la \gamma \ra$, where $A$ is isomorphic to a subgroup of $\Aut(\bbP_k^2)$. Since the latter has no elements of order 8, we are done.
\end{remark}

\section{Del Pezzo surfaces of degree 1}
This is the most difficult and interesting case. 
The linear system $|-2K_X|$ defines a degree 2 map $f:X\to Q$, where $Q$ is a quadric cone in $\bbP_k^3$. Again, since $-K_X$ is ample, $f$ is a finite  map, and arguing as in the previous case we see that the map is separable. The Galois group of the cover is generated by an automorphism $\beta$ of $X$ known as the \emph{Bertini involution}.  For any divisor $D$ we have 
\begin{equation}\label{eq2}
D+\gamma^*(D)\sim 2(D\cdot K_X)K_X.
\end{equation}
This shows that $\beta^*$ acts as the minus identity on the lattice $K_X^\perp$. The lattice $K_X^\perp$ is isomorphic to the root lattice of type $E_8$. The involution $\beta^*$ generates the center of the Weyl group $W(E_8)$.

The automorphism group $\Aut(X)$ is isomorphic to a subgroup of  $W(E_8)$. Possible orders $p^s, s > 1,$  of minimal automorphisms are $4$ and $8$ (see \cite{DI}). 

So we assume  $p = 2$. The  linear 
system $|-K_X|$ has one base point $p_0$.  Blowing it up we obtain a fibration $\pi:X'\to \bbP_k^1$ whose general fibre is an irreducible  curve of arithmetic genus 1. Since $-K_X$ is ample, all fibres are irreducible, and this implies that a general fibre is an elliptic curve (see \cite{CD}, Corollary 5.5.7). Let $S_0$ be the exceptional curve of the blow-up. It is a section of the elliptic fibration. We take  it as the zero in the Mordell-Weil group of sections of $\pi$.
The map $f:X\to Q$ extends to a degree 2 separable finite map $f':X'\to   \bfF_2$, where $\bfF_2$ is the minimal ruled surface with the exceptional section $E$ satisfying $E^2 = -2$. Its branch curve is equal to the union of $E$ and a curve $B$ from the divisor class $3f+e$, where $f$ is the class of a fibre and $e = [E]$. We have $f'{}^*(E) = 2S_0$.  The elliptic fibration on $X'$ is the pre-image of the ruling of $\bfF_2$. We know that  $\tau =\sigma^2$ acts identically on the base of the elliptic fibration. Since it also leaves invariant the section $S_0$, it defines an automorphism of the generic fibre considered as an abelian curve with zero section defined by $S_0$. If  $\tau^2 = 1$, then $\tau$ is the negation automorphism, hence defines the Bertini involution of $X$.  The group of automorphisms of an abelian curve in characteristic  $2$ is   of order 2 if the absolute invariant  of the curve is not equal to 0 or of order  24 otherwise. In the latter case it is isomorphic to $Q_8\rtimes \bbZ/3$, where $Q_8$ is the quaternion group with the center generated by the negation automorphism (see \cite{Silverman}, Appendix A). 
Thus $\tau^4 =  1$ and the Weierstrass model of the generic fibre is 
$$ y^2+a_3y+x^3+a_4x+a_6 = 0.$$
In global terms, the Weierstrass model of the elliptic fibration $\pi:X'\to \bbP_k^1$ is a surface  in $\bbP(1,1,2,3)$ given by the equation
$$y^2+a_3(u,v)y+x^3+a_4(u,v)x+a_6(u,v),$$
where $a_i$ are binary forms of degree $i$. It  is obtained by blowing down the section $S_0$ to the point $(u,v,x,y) = (0,0,1,1)$ and  is isomorphic to our Del Pezzo surface $X$. The image of the branch curve $B$ is given by the equation $a_3(u,v) = 0$, i.e. $B$ is equal to the pre-image of an effective divisor of degree 3 on the base plus the section $S_0$. Since a general point of $B$ is a 2-torsion point of a general fibre, we see that all nonsingular fibres of the elliptic fibration are supersingular elliptic curves (i.e. have no non-trivial 2-torsion points). An automorphism of order 4 of $X$ is defined by 
$$(u,v,x,y) \mapsto (u,v,x+ s(u,v)^2, y+s(u,v)x+t(u,v)),$$
where $s$ is binary forms of degree $1$ and $t$ is a binary form of degree 3 satisfying 
\begin{equation}\label{eq}
a_3 = s^3,\  t^2+a_3t+s^6+a_4s^2 = 0.
\end{equation}
In particular,  it shows that $a_3$ must be a cube, so we can change the coordinates $(u,v)$ to assume that $s = u, a_3 = u^3.$ The second equality in \eqref{eq} tells that $t$ is divisible by $u$, so we can write it as $t = uq$ for some binary form $q$ of degree 2 satisfying $q^2+u^2q+u^4+a_4 = 0$. Let $\alpha$ be a root of the equation $x^2+x+1=0$ and $b = q+\alpha u^2.$ Then $b$ satisfies
$a_4 = b^2+u^2b$ and $t = ub+\alpha u^3$. Conversely, any surface in $\bbP(1,1,2,3)$ with equation
\begin{equation}\label{serre}
y^2+u^3y+x^3+(b(u,v)^2+u^2b(u,v))x+a_6(u,v) = 0
\end{equation}
where $b$ is a quadratic form in $(u,v)$ and the coefficient at $uv^5$ in $a_6$ is not zero (this is equivalent to that the surface is nonsingular) is a Del Pezzo surface of degree 1 admitting an automorphism of order 4
$$\tau:(u,v,x,y) \mapsto (u,v, x+u^2,y+ux+ub+\alpha u^3).$$
Note that $\tau^2:(u,v,x,y) \mapsto (u,v,x,y+u^3)$ coincides with the Bertini involution.


\begin{theorem} Let $X$ be a Del Pezzo surface \eqref{serre}. Then it does not admit an automorphism of order 8.
\end{theorem}

\begin{proof} Assume $\tau = \sigma^2$. Since $\sigma$ leaves invariant $|-K_X|$, it fixes its unique base point, and lifts to an automorphism of the elliptic surface $X'$ preserving the zero section $S_0$.
Since the general fibre of the elliptic fibration $f:X' \to \bbP_k^1$ has no automorphism of order 8, the transformation $\sigma$ acts nontrivially on the base of the fibration. Note that the fibration has only one singular fibre $F_0$ over $(u,v) = (0,1)$. It is a cuspidal cubic. The transformation $\sigma$ leaves this fibre invariant and hence acts on $\bbP_k^1$ by $(u,v)\mapsto (u, u+cv)$ for some $c\in k$. Since the restriction of $\sigma$ to $F_0$ has at least two distinct  fixed points: the cusp and the origin $F_0\cap S_0$, it acts identically on $F_0$ and freely on its  complement $X'\setminus F_0$. 

Recall that $X'$ is obtained by blowing up 9 points  $p_1,\ldots,p_9$ in $\bbP_k^2$, the base points of a pencil of cubic curves. We may assume that $X$ is the blow-up of the first 8 points, and the exceptional curve over $p_9$ is the zero section $S_0$. Let $S$ be the exceptional curve over any other point.  We know that $\beta = \sigma^4$ is the Bertini involution of $X$. Applying formula \eqref{eq2}, we obtain that 
$S\cdot \beta(S) = 3.$ Identifying  $\beta(S)$  and $S$ with their pre-images in $X'$, we see that $\beta(S) + S  = S_0$ in the Mordell-Weil group of sections of 
$\pi:X'\to \bbP_k^1$. Thus $S$ and $\beta(S)$ meet at 2-torsion points of fibres. However, all nonsingular fibres of our fibration are supersingular elliptic curves, hence $S$ and $\beta(S)$ can meet only at the singular fibre $F_0$. Let $Q\in F_0$ be the intersection point. The sections $S$ and $\beta(S)$ are tangent to each other at $Q$ with multiplicity 3. Now consider the orbit of the pair $(S,\beta(S))$ under the cyclic group $\la \sigma\ra$. It consists of 4 pairs  
$$(S, \sigma^4(S)), \ (\sigma(S), \sigma^5(S)),\  (\sigma^2(S), \sigma^6(S)), \ (\sigma^3(S), \sigma^7(S)).$$
Let $D_i = \sigma^i(S)+\sigma^{i+4}(S), \ i = 1,2,3,4$. We have $D_1+\ldots+D_4 \sim -8K_X,$, hence for $i\ne j$ we have 
$D_i\cdot D_j = (64-16)/12 = 4$. Let $Y\to X$ be the blow-up of $Q$. Since $Q$ is a double point of each $D_i$, the proper transform $\bar{D}_i$ of each $D_i$ in $Y$ has self-intersection 0 and consists of two smooth rational curves intersecting at one point with multiplicity 2. Moreover, we have $\bar{D}_i\cdot \bar{D}_j = 0$.  Applying \eqref{eq2}, we get $D_i\in |-2K_X|$. Since $Q$ is a double point of $D_i$, we obtain $\bar{D}_i\in |-2K_{Y}|$. The linear system $|-2K_Y|$ defines a fibration $Y\to \bbP_k^1$ with a curve of arithmetic genus 1 as a general fibre (an elliptic or a quasi-elliptic fibration) and four singular fibres  $\bar{D}_i$  of Kodaira's type $III$. The automorphism $\sigma$ acts on the base of the fibration and the four special fibres form one orbit. But the action of $\sigma$ on $\bbP_k^1$ is of order 2 and this gives us a contradiction.
\end{proof}

\begin{remark} A computational proof of Theorem 7 was given by J.-P. Serre.
\end{remark}

\section{Conjugacy classes of elements of order $p^2$}

Assume that $k$ is algebraically closed. As we have seen in the previous sections, an element of order $p^2$ not conjugate to a projective transformation exists only for $p = 2$. It  can be realized as a minimal automorphism of a conic bundles,  or a  Del Pezzo surfaces of degree $1$ or $4$. Del Pezzo surfaces of degree $1$ are super-rigid, i.e. a minimal automorphism of such a surface could be conjugate only to a minimal automorphism of the same surface. A minimal automorphism of a Del Pezzo surface of degree 4 is conjugate to a minimal automorphism of a conic bundle with 5 singular fibres (see \cite{DI}, \S 8). 

Thus we have proved the following.

\begin{theorem}  An element of order $p^2$ not conjugate to a projective transformation exists only if $p = 2$.  Assume that $k$ is algebraically closed.  An element of order $4$ is  conjugate to either a projective transformation, or a transformation  realized by a minimal  automorphism of a conic  bundle or of a Del Pezzo surface of degree $1$.  
\end{theorem}

For the completeness sake let us add that elements of order $p$ not conjugate to a projective transformations occur for any $p$. They can be realized as automorphisms of conic bundles, and if $p = 2,3,5$ as automorphisms of Del Pezzo surfaces.

\bibliographystyle{amsplain}

\end{document}